\documentclass{article}

\usepackage{amsmath, amsthm, amssymb}

\newtheorem{theorem}{Theorem}[section]

\newtheorem{proposition}[theorem]{Proposition}

\begin{document}


\title{Interlacing of positive real zeros of Bessel functions}

\author{Tam\'as P\'almai\footnote{Electronic mail: palmai@phy.bme.hu} and Barnab\'as Apagyi\footnote{Electronic mail: apagyi@phy.bme.hu}\\
Department of Theoretical Physics\\Budapest University of Technology and Economics\\H-1111 Budapest, Hungary}

\maketitle

\begin{abstract}
We unify the three distinct inequality sequences (Abramowitz ans Stegun (1972) [1,9.5.2]) of positive real zeros of Bessel functions into a single
one.
\end{abstract}


\section{Introduction and results}

Consider the Bessel differential equation
\begin{equation}\label{1}
z^2\frac{d^2\mathcal{C}_\nu(z)}{dz^2}+z\frac{d\mathcal{C}_\nu(z)}{dz}+(z^2-\nu^2)\mathcal{C_{\nu}}(z)=0.
\end{equation}
The general solution is
$\mathcal{C}_\nu(z)=J_\nu(z)\cos\alpha-Y_\nu(z)\sin\alpha$, a
linear combination of $J_\nu(z)$ and $Y_\nu(z)$ being the Bessel
functions of the first and second kind (defined e.g. in Ref.
\cite{Abramowitz}).

When $\nu$ is real, the Bessel functions $J_{\nu}(x)$,
$Y_{\nu}(x)$, their derivatives $J'_{\nu}(x)$, $Y'_{\nu}(x)$ each
have an infinite number of positive real zeros, all of which are
simple with the possible exception of $x=0$ \cite{Abramowitz} (see
also \cite{Watson}, Bessel-Lommel and Rolle's theorems). For
non-negative order $\nu$ the $s$th positive real zeros of these
functions are denoted by $j_{\nu,s}$, $y_{\nu,s}$, $j'_{\nu,s}$,
$y'_{\nu,s}$, respectively, except that $x=0$ is counted as the
first zero of $J'_{0}(x)$. Since $J'_{0}(x)=-J_{1}(x)$, it follows
that $j'_{0,1}=0,j'_{0,s}=j_{1,s-1}, (s=2,3,...)$
\cite{Abramowitz}.

The following results are widely known
\cite{Abramowitz,Watson,LiuZou} about the zeros of the Bessel
functions.
\begin{theorem}\label{thm:1}
The positive real zeros of the functions $J_\nu(x)$,
$J_{\nu+1}(x)$, $Y_\nu(x)$, $Y_{\nu+1}(x)$, $J'_\nu(x)$ and
$Y'_\nu(x)$ interlace according to the three distinct
inequalities:
\begin{equation}\label{3}
j_{\nu,1}<j_{\nu+1,1}<j_{\nu,2}<j_{\nu+1,2}<j_{\nu,3}<\ldots
\end{equation}
\begin{equation}\label{4}
y_{\nu,1}<y_{\nu+1,1}<y_{\nu,2}<y_{\nu+1,2}<y_{\nu,3}<\ldots
\end{equation}
\begin{equation}\label{5}
\nu\leq j'_{\nu,1}<y_{\nu,1}<y'_{\nu,1}<j_{\nu,1}<j'_{\nu,2}
<y_{\nu,2}<y'_{\nu,2}<j_{\nu,2}<j'_{\nu,3}<\ldots
\end{equation}
Furthermore, the positive real zeros of $J'_\nu(x)$, $J'_{\nu+1}(x)$ and that of $Y'_\nu(x)$, $Y'_{\nu+1}(x)$ are interlaced:
\begin{equation}\label{lz1}
j'_{\nu,1}<j'_{\nu+1,1}<j'_{\nu,2}<j'_{\nu+1,2}<j'_{\nu,3}<\ldots
\end{equation}
\begin{equation}\label{lz2}
y'_{\nu,1}<y'_{\nu+1,1}<y'_{\nu,2}<y'_{\nu+1,2}<y'_{\nu,3}<\ldots
\end{equation}
\end{theorem}
Note that Eqs. (\ref{lz1})-(\ref{lz2}) were established only
recently in Ref. \cite{LiuZou} where a particular inverse
scattering problem was studied.

In addition to Theorem \ref{thm:1} the following auxiliary relations can be found.

\begin{proposition}\label{prop}
For non-negative orders, i.e. $\nu\geq0$
\begin{eqnarray}
j_{\nu+1,s}<j'_{\nu,s+1},\quad s=1,2,\ldots\label{palz}\\
y_{\nu+1,s}<y'_{\nu,s},\quad s=1,2,\ldots\label{propeq}
\end{eqnarray}
hold.
\end{proposition}
Eq. (\ref{palz}) emerged previously when studying a particular inverse scattering problem \cite{PA3}, and also in Ref. \cite{LiuZouTR} independently of the authors.

Now, with the aid of Proposition \ref{prop} it is possible to unify Eqs. (\ref{3}), (\ref{4}) and (\ref{5}) into a single one. In addition we obtain a simple proof of Eqs. (\ref{lz1}) and (\ref{lz2}).

We shall formulate our main result in a slightly generalized way, including also an interesting breaking condition.
\begin{theorem}[Interlacing of positive real zeros of the Bessel functions]\label{thm:3}
The positive real zeros of the Bessel functions $J_{\nu}(x)$,
$J'_{\nu}(x)$, $Y_{\nu}(x)$, $Y'_{\nu}(x)$,
$J_{\nu+\varepsilon}(x)$, $Y_{\nu+\varepsilon}(x)$,
$0<\varepsilon\leq1$, are interlaced according to the inequalities
\begin{equation}\label{ineq}
j'_{\nu,s}<y_{\nu,s}<y_{\nu+\varepsilon,s}<y'_{\nu,s}<j_{\nu,s}<j_{\nu+\varepsilon,s}<j'_{\nu,s+1}
<\ldots\quad s=1,2,\ldots,\,\nu\geq 0.
\end{equation}
For $\varepsilon>1$ this interlacing property is destroyed.
\end{theorem}

Eqs. (\ref{lz1}) and (\ref{lz2}) can be generalized to
\begin{equation}\label{jder}
j'_{\nu,1}<j'_{\nu+\varepsilon,1}<j'_{\nu,2}<j'_{\nu+\varepsilon,2}<j'_{\nu,3}<\ldots
\end{equation}
\begin{equation}\label{yder}
y'_{\nu,1}<y'_{\nu+\varepsilon,1}<y'_{\nu,2}<y'_{\nu+\varepsilon,2}<y'_{\nu,3}<\ldots
\end{equation}
with $0<\varepsilon\leq1$ and $\nu\geq0$. We note that the latter
two inequalities cannot be integrated with our unified interlacing
inequality (\ref{ineq}). While both
$j'_{\nu,s}<j'_{\nu+\varepsilon,s}<y_{\nu+\varepsilon,s}$ and
$y'_{\nu,s}<y'_{\nu+\varepsilon,s}<j_{\nu+\varepsilon,s}$ hold,
numerical counterexamples can easily be constructed for the
non-existence of a uniform inequality between both
$j'_{\nu+\varepsilon,s}$ and $y_{\nu,s}$ (for which
$j'_{\nu,s}<y_{\nu,s}<y_{\nu+\varepsilon,s}$ applies), and
$y'_{\nu+\varepsilon,s}$ and $j_{\nu,s}$ (for which
$y'_{\nu,s}<j_{\nu,s}<j_{\nu+\varepsilon,s}$ applies) for all
$s=1,2,\ldots$ and $\nu>0$, $0<\varepsilon\leq1$.

\section{Proofs}

We start by proving Eqs. (\ref{jder}) and (\ref{yder}) depending on Theorem \ref{thm:3}. The
first inequality is trivial due to the monotonicity of $j'_{\nu,s}$
in $\nu$ for $\nu\geq0$ (see \cite{Watson}). For proving the second inequality take the sequence of
Theorem \ref{thm:3} at some arbitrary $s$ for $\nu$, $\nu+1$ and for
$\nu+1$, $\nu+2$ with $\varepsilon=1$:
\begin{eqnarray*}
j'_{\nu,s}<y_{\nu,s}<y_{\nu+1,s}<y'_{\nu,s}<j_{\nu,s}<j_{\nu+1,s}<j'_{\nu,s+1}\\
j'_{\nu+1,s}<y_{\nu+1,s}<y_{\nu+2,s}<y'_{\nu+1,s}<j_{\nu+1,s}<j_{\nu+2,s}<j'_{\nu+1,s+1}
\end{eqnarray*}
From the first one we have $j_{\nu+1,s}<j'_{\nu,s+1}$ and from the
second one we have $j'_{\nu+1,s}<j_{\nu+1,s}$. Combining these two
yields the first interlacing property $j'_{\nu+1,s}<j'_{\nu,s+1}$
of Eq. (\ref{lz1}) for the derivative function $J'_\nu(x)$. For
$\varepsilon<1$ the relation follows from the the monotonicity of
$j'_{\nu,s}$ in $\nu$.

For the positive zeros of the derivative function $Y'_\nu(x)$ a
similar reasoning can be presented.
$y'_{\nu,s}<y'_{\nu+\varepsilon,s}$ is trivial due to the
monotonicity in $\nu$. Use Theorem \ref{thm:3} for $\nu$, $\nu+1$ and for
$\nu+1$, $\nu+2$ with $\varepsilon=1$
\begin{eqnarray*}
\ldots<y_{\nu+1,s+1}<y'_{\nu,s+1}<\ldots\\
\ldots<y'_{\nu+1,s}<j_{\nu+1,s}<j'_{\nu+1,s+1}<y_{\nu+1,s+1}<\ldots
\end{eqnarray*}
and combine them to get the relations $y'_{\nu+1,s}<y'_{\nu,s+1}$
contained in Eq. (\ref{lz2}) . Again the monotonicity of
$y_{\nu,s}$ in $\nu$ implies the non-trivial inequalities
$y'_{\nu+\varepsilon,s}<y'_{\nu,s+1}$ for $\varepsilon<1$.

Eq. (\ref{palz}) of Proposition  \ref{prop} was already proven,
independently of each other, in Refs. \cite{PA3} and
\cite{LiuZouTR}  therefore its proof is omitted here. The proof of
Eq. (\ref{propeq}) is elementary and based on the analysis of
intervals on which both $Y_{\nu+1}(x)$ and $Y_\nu(x)$ take the
same and the opposite sign. Of course, by definition,
$Y_{\nu+1}(x)$ and $Y_{\nu}(x)$  each keeps sign in the intervals
defined by their two consecutive zeros, respectively. That is
$Y_{\nu+1}(x)$ keeps the sign in the interval
\begin{equation*}
y_{\nu+1,s}<x<y_{\nu+1,s+1},\qquad s=1,2,\ldots,
\end{equation*}
and $Y_{\nu}(x)$ does it in
\begin{equation*}
y_{\nu,s}<x<y_{\nu,s+1},\qquad s=1,2,\ldots.
\end{equation*}
However, since $Y_{\nu}(x\to 0)=-\infty$ for $\nu\geq0$ and using
Eq. (\ref{4}) [which implies that $y_{\nu,s}<y_{\nu+1,s} $ and
$y_{\nu,s+1}<y_{\nu+1,s+1} $], one concludes by induction that
both $Y_{\nu+1}(x)$ and $Y_\nu(x)$ keep the same sign in the
common intervals
\begin{equation}\label{17}
y_{\nu+1,s}<x<y_{\nu,s+1},\qquad s=1,2,\ldots,
\end{equation}
whereas in
\begin{equation}\label{18}
y_{\nu,s}<x<y_{\nu+1,s},\qquad s=1,2,\ldots,
\end{equation}
the signs do differ. Now let us take the recurrence relation
$$
\mathcal{C}'_\nu(x)=-\mathcal{C}_{\nu+1}(x)+\frac{\nu}{x}\mathcal{C}_\nu(x)
$$
with $\mathcal{C}=Y$ at $x=y'_{\nu,s}$. It yields
\begin{equation}\label{19}
Y_{\nu+1}(y'_{\nu,s})=\frac{\nu}{y'_{\nu,s}}Y_{\nu}(y'_{\nu,s}),
\end{equation}
i.e. the signs of  $Y_{\nu+1}(x)$ and $Y_{\nu}(x)$ coincide at
$x=y'_{\nu,s}$. But, because of (\ref{5}) [which tells that
$y'_{\nu,s}$ lies within the interval $y_{\nu,s}<x<y_{\nu,s+1}$],
the content of Eq. (\ref{19}) means also that $y'_{\nu,s}$ must be
in the common intervals given above by  (\ref{17}). This completes
the proof of Proposition \ref{prop} for $y_{\nu+1,s}<y'_{\nu,s}$.

The proof of the first part of Theorem \ref{thm:3} is also elementary and follows from the
application of the two relations of Proposition \ref{prop} [being
previously unknown] in conjunction with the three distinct
inequalities (\ref{3}), (\ref{4}) and (\ref{5}) [being  already known, i.e. from
Ref. \cite{Abramowitz}]. The case $0<\varepsilon<1$ is  immediately implied
by the  well-known property of $j_{\nu,s}$'s, and $y_{\nu,s}$'s that, for a
fixed $s$, they are strictly increasing functions of $\nu$ if $\nu\geq0$ \cite{Watson}.

 The negative statement for $\varepsilon>1$ can be deduced from
$y_{\nu+\varepsilon,s}>j_{\nu,s}$, that is from the violation of
the prescribed relation between the third and fifth term in the
inequality sequence of Theorem \ref{thm:3}. In Ref. \cite{PA3} it
was proven that the Wronskian  $W_{\nu\mu}(x)\equiv
J_\nu(x)Y_\mu'(x)-J_\nu'(x)Y_\mu
(x)\neq0$ for $x\in(0,\infty)$ if
and only if $0<|\nu-\mu|\leq1$ is maintained ($\nu\neq\mu$). One
of the ideas in that proof is that the set of extremal points of
$W_{\nu\mu}(x)$ is
$\{j_{\nu,s}\}_{s=1}^\infty\cup\{y_{\mu,s}\}_{s=1}^\infty$ and it
has been unveiled that the inequality sequences
$y_{\mu,s}<j_{\nu,s}<y_{\mu,s+1}<j_{\nu,s+1}$, $s=1,2,\ldots$
 hold
if and only if all the local extrema are of the same sign, i.e.
$W_{\nu\mu}(x)\neq0$, $x\in\mathbb{R}^+$.
These relations are exactly the same that we are studying here
with $\mu=\nu+\varepsilon$. Since for $\varepsilon>1$
$|\nu-\mu|\leq1$ cannot hold and there is at least one root of
$W_{\nu\mu}(x)$ the inequalities must be violated for some $s$.
Now our proof is complete.


\begin{thebibliography}{99}

\bibitem{Abramowitz}
Abramowitz M and Stegun I A: {\it Handbook of Mathematical
Functions} (New York: Dover Publications) pp. 360-371 (1972)

\bibitem{PA3}
P\'almai T and Apagyi B: On nonsingular potentials of Cox-Thompson
inversion scheme, {\it Journal of Mathematical Physics}
{\bf 51}, 022114 (2010)


\bibitem{Watson}
Watson G N: {\it A treatise on the theory of Bessel functions}
(Cambridge Mathematical Library edition) chapter 15 (1995)

\bibitem{LiuZou}
Liu H Y and Zou J: Zeros of the Bessel and spherical Bessel
functions and their applications for uniqueness in inverse
acoustic obstacle scattering, {\it IMA Journal of Applied Mathematical}
{\bf 72}, pp. 817-831 (2007)

\bibitem{LiuZouTR}
Liu H Y and Zou J: Zeros of the Bessel and spherical Bessel
functions and their applications for uniqueness in inverse
acoustic obstacle scattering problems, {\it Technical Report CUHK-2007-02 (342)},
The Chinese University of Hong Kong, Hong Kong, 2007.


\end{thebibliography}
\end{document}